\documentclass[11pt]{article}
\usepackage{latexsym}

\newtheorem{theorem}{Theorem}{}
{}
{}

\newcommand{\R}{\mathcal{R}}

\begin{document}

\title{Interval Semantics\\for Standard Floating-Point Arithmetic}

\date{}

\author{W.W. Edmonson \\ North Carolina State University
        \and
        M.H. van Emden \\ University of Victoria
}
\date{\small{
        Research Report DCS-323-IR \\
        Department of Computer Science, University of Victoria
      }
}

\maketitle

\begin{abstract}
If the non-zero finite floating-point numbers are interpreted as point
intervals, then the effect of rounding can be interpreted as computing
one of the bounds of the result according to interval arithmetic.
We give an interval interpretation for the signed zeros and infinities,
so that the undefined operations $\pm 0 * \pm \infty$,
$\pm \infty - \pm \infty$, 
$\pm \infty / \pm \infty$, 
and
$\pm 0 / \pm 0$
become defined.

In this way no operation remains that gives rise to an error
condition.
Mathematically questionable features of the floating-point standard
become well-defined sets of reals.
Interval semantics provides a basis for the verification of
numerical algorithms.
We derive the results of the newly defined operations
and consider the implications for hardware implementation.

\end{abstract}


\section{Introduction}
IEEE Standard 754-1985 for Binary Floating-Point Arithmetic
\cite{Std754}\footnote{Its successor, IEEE Standard 754-2008 for
  Floating-Point Arithmetic\cite{Std754R}, does not make changes that are
  relevant to this paper.} achieved great success by finding a
synthesis of the best features of the existing processors and causing
these to be widely adopted in a short time.  It was not to be expected
that IEEE Standard 754 had a remedy for the fact that none of the
existing processors was based on a coherent approach to the
fundamental problem of how to approximate on digital processors the
operations of arithmetic on reals.  An example of the ad-hoc approach
taken in IEEE Standard 754 is the introduction of the infinities.
Mathematically, there is one advantage: division of a non-zero number
by zero becomes defined.  But it introduces more exceptions than it
removes: $0 * \pm \infty$, $\pm \infty - \pm \infty$, $\pm \infty /
\pm \infty$, while $0/0$ remains undefined.

Any improvement to the standard needs to be based on a mathematically
convincing approach to the following problem:
\begin{quote}
\emph{
How to map the
arithmetical structure of the reals to a closed, exception-free algebra
on a set of computer-representable quantities.
}
\end{quote}
In this paper we review results from \cite{hckvnmdn01} (see there for
earlier references) that allow such an algebra to be based on
intervals.  Surprisingly, most of the standard carries over unchanged
to our algebra. What does not carry over are the error conditions:
there are none.

Our proposal is based on the idea of interpreting floating-point
numbers as \emph{sets of reals}.  The sets of reals include the reals
themselves by identifying the singleton sets $\{x\}$ with $x$ itself,
for all reals $x$.  As we will show, this defines an arithmetic that
is compatible with arithmetic on the reals.  Where it deviates is that
division by zero becomes defined.

But there is a more important advantage.
Any arithmetic that intends to approximate real arithmetic with a
finite set of values necessarily introduces \emph{uncertainty}.
For example, in binary floating-point arithmetic the result of dividing 1 by 10
leaves uncertainty concerning the digits from a certain point onwards.
Also, as there is necessarily a greatest floating-point number $M$,
there has to be uncertainty about any result that exceeds this
maximum.
Hence there need to be floating-point numbers that are interpreted as
$\{ x \in \R \mid M \leq x \}$
and as
$\{ x \in \R \mid -M \geq x \}$.

In Section~\ref{preliminaries} we show how the theoretical advantages
of arithmetic on sets of reals becomes practical by restricting these
sets to be intervals. We also establish notation and terminology for
intervals and floating-point numbers.  In Section~\ref{semFlpt} we
define the central feature of this paper: an interval interpretation
of each floating-point number.  Here we take the view that zero has
finite precision.  In Section~\ref{fragmIA} we state theorems that
establish that to a large extent the existing floating-point standard
already implements the arithmetic operations on floating-point numbers
interpreted as intervals.  Section~\ref{mods} treats the case where
interval semantics changes the existing definition;
Section~\ref{undefDef} derives the undefined cases according to
interval semantics.  Section~\ref{altSem} is a short discussion of
infinite-precision zeroes.  Finally, in Section~\ref{hardwCons} we
survey the consequences of the new semantics for hardware
implementation.

\section{Preliminaries}
\label{preliminaries}

\subsection{Floating-point numbers}
If $x$ is a floating-point number greater than $-\infty$ (less than
$+\infty$), then $x^-$ ($x^+$) is the next smaller (greater)
floating-point number. The constant $m$ is defined as the least
positive floating-point number and $M$ as the greatest finite
floating-point number.

\subsection{Intervals}

``Real intervals'' are to be interpreted according to the following
definition: a \emph{real interval} is a \emph{closed, connected set of
  reals}.
According to a well-known result in topology,
such sets take the following forms:
$\{x \in \R \mid x \leq b \}$,
$\{x \in \R \mid a \leq x \leq b \}$,
$\{x \in \R \mid a \leq x \}$,
$\R$ and $\emptyset$.
Here $a$ and $b$ are reals such that $a \leq b$.
We denote these sets as $(-\infty, b]$,
$[a,b]$, $[a, +\infty)$, $(-\infty, +\infty) = \R$, and $\emptyset$.

If $x$ is an interval, $x_l$, the \emph{left bound} of $x$,
is the greatest lower bound of $x$ as a set of reals,
if it has one, otherwise $x_l$ is $-\infty$;
$x_r$, the \emph{right bound} of $x$,
is the least upper bound of $x$ as a set of reals,
if it has one, otherwise $x_r$ is $+\infty$.
An interval may consist of a single real number $s$, hence written
$[s,s]$.  This is called a \emph{point interval}, sometimes also
called \emph{degenerate interval}.

\paragraph{Arithmetic on real intervals}

The na\"{i}ve approach to interval arithmetic is embodied in the
following tentative definition:
for all intervals $X$ and $Y$,
\begin{eqnarray*}
X + Y & = & \{ x+y \in \R \mid x \in X \wedge y \in Y \}   \\
X - Y & = & \{ x-y \in \R \mid x \in X \wedge y \in Y \}   \\
X * Y & = & \{ x*y \in \R \mid x \in X \wedge y \in Y \}   \\
X / Y & = & \{ x/y \in \R \mid x \in X \wedge y \in Y \}   \\
\end{eqnarray*}
The na\"{i}ve approach gets everything right
except for division when $0 \in Y$.

This problem can be remedied by invoking in the interval arithmetic
operations, the definitions of the inverse operation $/$ on the reals
as a multiplication.  Even though this is not necessary for the other
operations, in the following definition (still tentative) these
operations have been modified in the same way.
\begin{eqnarray*}
X + Y & = & \{ z \in \R \mid \exists x \in X,  y \in Y \;.\; x+y = z \}   \\
X - Y & = & \{ z \in \R \mid \exists x \in X,  y \in Y \;.\; y+z = x \}   \\
X * Y & = & \{ z \in \R \mid \exists x \in X,  y \in Y \;.\; x*y = z \}   \\
X / Y & = & \{ z \in \R \mid \exists x \in X,  y \in Y \;.\; y*z = x \}   \\
\end{eqnarray*}
The operations defined in this way have the advantage of being defined
for all intervals.
But the division operation again gives a problem:
although now division is everywhere defined,
it may fail to yield an interval, as in
$[1,1] / [-1,+1] = (-\infty, -1] \cup [+1, +\infty)$.
This is not an acceptable result, as we need a \emph{closed},
exception-free algebra.

As we insist on the simplicity of intervals for the inclusion of the
sets of reals that arise in computation,
we need to enclose some of these sets in the smallest
interval containing them.
Hence the introduction of the \emph{interval hull} operator $\Box$
where $\Box S$ is defined as the least interval containing $S$, for any
set $S$ of reals.

For the sake of uniformity we include the interval hull operator also
in those cases where its argument is already an interval. Hence the
following, definitive, definition:
\begin{eqnarray}
\label{intvOps}
X + Y & = &
\Box \{ z \in \R \mid \exists x \in X,  y \in Y \;.\; x+y = z \}
\nonumber \\
X - Y & = &
\Box \{ z \in \R \mid \exists x \in X,  y \in Y \;.\; y+z = x \} \\
X * Y & = &
\Box \{ z \in \R \mid \exists x \in X,  y \in Y \;.\; x*y = z \}
\nonumber \\
X / Y & = &
\Box \{ z \in \R \mid \exists x \in X,  y \in Y \;.\; y*z = x \}
\nonumber
\end{eqnarray}
In the general case it is quite complicated, for the multiplication
and division formulas, to express the bounds of the result intervals
in terms of operations on the operand bounds. In this paper we only
need to consider the case where at least one operand is a point
interval. Another simplifying circumstance is that, when the bounds of
an operand are not equal, their signs are equal.

With definition (\ref{intvOps}), real interval arithmetic is a closed,
exception-free algebra. It is defined for all interval operands and
always yields an interval.

\paragraph{Floating-point intervals}

A floating-point interval is a real interval where the bounds are
restricted to floating-point numbers. Every finite set $F$ of reals
defines a finitary approximation to real arithmetic by mapping every
real interval to the least floating-point interval containing it.

The interval hull introduced above is with respect to the set of real
intervals. From now on we are only concerned with the set of
floating-point intervals. Accordingly, the interval hull $\Box$ in
Equations (\ref{intvOps}) will be taken with respect to the
floating-point intervals. The arithmetic operations define a closed,
exception-free algebra on floating-point intervals.

\section{A proposed semantics for floating-point numbers}
\label{semFlpt}

We propose to interpret floating-point numbers as intervals according
to the following cases.
\begin{itemize}

\item Each non-zero finite floating-point number $x$ is interpreted as
  the point interval $[x,x]$; that is, the singleton set
  $\{x\}$.

\item The floating-point number $+\infty$ is interpreted as the real
  interval $[M,+\infty)$.

  As $M$ is the greatest finite floating-point number, there is
  uncertainty about any result greater than that.  Hence, it is
  desirable to have a notation for the set of all those values. For
  technical reasons $M$ is included in this set\footnote{
This is mathematically correct (the result is included in the set),
but not as informative as when the lower bound of the interval were
allowed to be open. This we don't allow so as to keep interval
representations simple.}.
The floating-point
  number $+\infty$, being different from any finite floating-point
  number, is an appropriate representation for $[M,+\infty)$.

\item The floating-point number $-\infty$ is interpreted as
  $(-\infty,-M]$ for reasons similar to the previous case.

\item The floating-point number $-0$ is interpreted as the
  floating-point interval $[-m,0]$, where $m$ is the least positive
  floating-point number.  The floating-point number $+0$ is
  interpreted as the floating-point interval $[0,m]$.

  It might be thought that $[0,0]$ is the most appropriate
  interpretation. Consider, however, the result of $a/0$, which would
  be interpreted as $[a,a]/[0,0]$, and which is the empty set for
  finite non-zero $a$ according to Equation (\ref{intvOps}).  The
  interpretation proposed here follows the idea that zero is a kind of
  dual to infinity: infinitely small versus infinitely large. Just as
  we do not try to give infinite precision to the infinitely large,
  this interpretation of the zeroes renounces infinite precision for
  the infinitely small. Thus we refer to this semantics as
  \emph{finite-precision zeroes}.  In Section~\ref{mods} it will be
  seen that this interpretation of the zeros interacts well with the
  other interval interpretations.  In Section~\ref{altSem} we explore
  the ramifications of infinite-precision zeroes.
\end{itemize}

\section{IEEE Standard 754 already implements a fragment of interval
  arithmetic}
\label{fragmIA}

The justification of our proposal is that for finite floating-point operands,
IEEE Standard 754 already implements interval arithmetic, in the sense of the
following theorem.
\begin{theorem}
\label{mainThm}
Suppose that the rounding mode is towards $+\infty$
$(-\infty)$ and that $x$ and $y$ are finite floating-point numbers.
The result of the floating-point operation $x \circ y$
is the upper (lower) bound of $[x,x] \circ [y,y]$
as computed in interval arithmetic, where $\circ$ is any of $+$, $-$,
and $*$.
\end{theorem}
\emph{Proof}:
The effect of the rounding mode is described as follows in the
IEEE Standard 754 \cite{Std754}:
\begin{quote}
When rounding toward +INFINITY the result shall be the format's value
(possibly +INFINITY) closest to and no less than the infinitely
precise result.
\end{quote}
With slightly different wording IEEE Standard 754R says \cite{Std754R}:
\begin{quote}
With $\ldots$ roundTowardPositive, the result shall be the format's
floating-point number (possibly $+\infty$) closest to and no less than
the infinitely precise result
\end{quote}
The ``infinitely precise result'' is the result according to real
arithmetic. The ``format's value closest to and no less than'' is the
right bound of the interval resulting from applying the interval hull
operator when the ``format's values'' constitute the floating-point
numbers used in the approximation operator.  A similar reasoning
applies when the rounding mode is toward $-\infty$.  $\Box$

A similar result holds when the rounding mode is toward $0$.
No such result holds when rounding is towards nearest,
as it is not known in which direction rounding takes place.
This issue is addressed in Section~\ref{hardwCons}.

\begin{theorem}
\label{auxThm}
Suppose that the rounding mode is towards $+\infty$
($-\infty$) and that $x$ and $y$ are finite floating-point numbers
with $y \neq 0$.
The result of the floating-point operation $x / y$
is the upper (lower) bound of $[x,x] / [y,y]$
as computed in interval arithmetic.
\end{theorem}

The statement of this theorem and its proof are similar to those of
Theorem~\ref{mainThm}. The present theorem is only separate because of
the need to exclude the possibility that $y = 0$.

These two theorems express what we mean when we say that IEEE Standard
754 already implements a fragment of interval arithmetic.

Interesting things happen when we allow $x$ or $y$ to become infinite,
or $y$ to become $0$ in $x/y$. These conditions will be discussed in
Section~\ref{undefDef}.

\section{Modification of interval arithmetic resulting from interval
  semantics}\label{mods}

Assigning, as we have done, an interval as meaning to every floating-point
number has to cause deviations from IEEE Standard 754.
This is because in interval arithmetic every operation has a defined
result, and this is not the case in IEEE Standard 754.
Under interval semantics some operations that are defined change;
some undefined ones become defined.
In this section we review the operations that are defined in IEEE
Standard 754.

\begin{itemize}
\item
$+\infty * +\infty = [M,+\infty)*[M,+\infty) = \Box [M^2, +\infty)$
according to floating-point interval arithmetic.
As $M^2$ is greater than the greatest floating-point number,
the $\Box$ operator widens its argument so that we get
$+\infty * +\infty = [M, +\infty)$.
The resulting interval is represented by the floating-point number $+\infty$.
Thus interval semantics maintains the result according to IEEE Standard 754.
\item
$+\infty * -\infty = (-\infty, -M] \sim -\infty$ according to similar
reasoning as above.
\item
For a positive finite floating-point number $a$,
$a * +\infty = [a,a] * [M,+\infty)$.
This evaluates to $[M,+\infty) \sim +\infty$ if $a \geq 1$,
otherwise to $[(a*M)^-, +\infty)$.
Thus under interval semantics $a * +\infty$ evaluates, if $a < 1$, to
$(a*M)^-$ or to $+\infty$, depending on rounding mode.
The IEEE Standard 754 result is maintained by interval semantics if $a \geq 1$.
It is not maintained if $a$ is sufficiently less than $1$.
We see that IEEE Standard 754 acts as if $\infty$ is infinitely precise.
According to interval semantics only selected finite reals are
infinitely precise, and $\infty$ has limited precision,
and this shows by multiplying it with a sufficiently small number.
\item
$\infty + \infty = [M,\infty) + [M,+\infty) = \Box [M+M,\infty)
=[M,\infty) \sim \infty$,
which conforms to IEEE Standard 754.
\item
For finite floating-point number $a$,
$a+\infty = [a,a]+[M,\infty) = [min((a+M)^-, M),+\infty)$.
Thus it conforms to IEEE Standard 754 for nonnegative $a$.
As $a$ approaches $-M$, the  result approaches $[0,+\infty)$.
\item
$+0 + (+0) = [0,m] + [0,m] = [0,2m]$
\item
$+0 + (-0) = [0,m] + [-m,0] = [-m,m]$
\item
For finite positive floating-point $a$,
$$a/+\infty = [a,a]/[M,+\infty) = \Box (0,(a/M)^+] = [0,(a/M)^+].$$
For moderate values of $a$, this is close to the proposed interval
interpretation of $+0$. For very large $a$, it approaches to $[0,1]$.
This is because $+\infty$ harbours a considerable degree of
uncertainty, as it stands for all reals greater than $M$.
\item
For finite positive floating-point $a$,
$$+\infty/a = [M,+\infty)/[a,a] = [min(M, (M/a)^-), \infty).$$
Departs moderately from IEEE Standard 754, and does so in a meaningful
way.  For $a$ equal to one, it is the interval interpretation of
$+\infty$.  As $a$ increases, the lower bound of the result of
$\infty/a$ decreases until the result becomes $[1,+\infty)$ for the
greatest finite value of $a$.  This is in accordance with the notion
that the interval interpretation of $\infty$ is not infinitely
precise.
\item $+\infty/+0 = [M,+\infty)/[0,m] = \Box [(M/m),+\infty) =
  [M,+\infty) \sim +\infty$,
which gives the IEEE Standard 754 result by mathematical reasoning rather than
by fiat.
\item
For positive finite $a$,
$a/+0 = [a,a]/[0,m] = [(a/m)^-, \infty)$.
For $a$ not very small, this result is about the interval
interpretation of $+\infty$, which is close to what IEEE Standard 754 defines.
As $a$ approaches $m$, this approaches $[1,+\infty)$.
\end{itemize}

\section{Undefined operations become defined}
\label{undefDef}

In Theorems~\ref{mainThm} and \ref{auxThm} we showed that for ordinary
cases interval semantics conform to IEEE Standard 754 outcomes.  In
Section~\ref{mods} we reviewed cases involving infinity where IEEE
Standard 754 does specify a result, but where that result may be
different under interval semantics. It remains to review the cases
where IEEE Standard 754 specifies that the result is NaN. As interval
arithmetic is everywhere defined, interval semantics specifies a
result in these cases.
\begin{itemize}
\item
$0 * +\infty = [0,m]*[M,+\infty) = [0,+\infty)$
\item
$+\infty/+\infty = [M,+\infty)/[M,+\infty) = \Box((0,+\infty)) = [0,+\infty)$
\item
$+\infty/-\infty = [M,+\infty)/(-\infty,-M] = \Box((-\infty,0)) = (-\infty,0]$
\item
$-\infty/-\infty = (-\infty,-M]/(-\infty,-M] = \Box((0,+\infty)) =
[0,+\infty)$
\item
$+0/+0 = \Box(\{z \in \R  \mid \exists x \in [0,m], y \in [0,m] \;.\;
yz = x \}) = [0,+\infty)$
\item
$+\infty - +\infty = [M,+\infty) - [M,+\infty) = (-\infty, +\infty) = \R$
\end{itemize}

\section{Infinite-precision zeroes}
\label{altSem}

In this section we investigate the consequences of interpreting the
floating-point numbers $+0$ and $-0$ both as the interval $[0,0]$.  An
important difference is that for positive finite floating-point $a$,
$a/+0$ becomes $[a,a]/[0,0]$ which equals
$$\Box(\{z \in \R  \mid \exists x \in [a,a], y \in [0,0] \;.\;
yz = x \}) = \emptyset.$$
With finite-precision zeroes the empty interval never arises.

With either kind of zeroes, every operation has a defined outcome, and
we have no need for NaN.  But with infinite precision zeroes, we
now have a likely candidate for an interval semantics for NaN: the
empty set. Strong support for this idea comes from the fact that,
whenever one of the operands is an empty interval, the result is
empty. That is, the empty interval propagates in the same way as
NaN. Also, whenever an empty interval arises, the programmer is likely to
want to treat this as an exceptional condition (though it is not an
error). Thus, under interval semantics, there is no INVALID OPERATION,
as there is in IEEE Standard floating-point arithmetic. Yet this flag
arises in the same situations, but now means the exceptional condition
of an empty interval being the result.

In the remainder of this section we list the formulas involving zero of
Sections~\ref{mods} and \ref{undefDef} under infinite-precision zeroes.

From Section~\ref{mods}:
\begin{itemize}
\item
$+0 + (+0) = [0,0] + [0,0] = [0,0]$
\item
$+0 + (-0) = [0,0] + [0,0] = [0,0]$
\item $+\infty/+0 = \Box(\{z \in \R \mid \exists x \in [M,+\infty), y
  \in [0,0] \;.\; yz = x \}) = \emptyset $
\item
For positive finite $a$,
$a/+0 = \emptyset$, as derived above.
\end{itemize}

From Section~\ref{undefDef}:
 
\begin{itemize}
\item
$0 * +\infty = [0,0]*[M,+\infty) = [0,0]$.
\item
$+0/+0 = +0/-0 =
\Box(\{z \in \R  \mid \exists x \in [0,0], y \in [0,0] \;.\;
yz = x \}) = (-\infty,+\infty)$
\end{itemize}

\section{Consequences for hardware implementation}
\label{hardwCons}

If we choose for the zeroes the finite precision semantics, there is
no role for the NaNs.  It might seem an advantage for hardware
implementation that NaN, its associated logic, exception, and flag can
be omitted. However, there seems to be no scarcity in fully IEEE
Standard 754 compliant hardware implementations. Indeed, most of the
logic is needed for the implementation of the algorithms for the
operations, so that exploiting the absence of NaNs provides only minor
relief.

An intriguing possibility is to do full interval arithmetic in
hardware \cite{iahpsun,hwsupia,vpia}, but that is not the topic of this paper.
General interval arithmetic is complex because the bounds of the same
interval may differ in sign.
It depends on their signs which bounds of the interval operands
should be used.  In this paper, we are concerned
with the simple case of operating on point intervals so that the
resulting interval is a point interval or it has adjacent floating-point
numbers as bounds.  As a result, if the upper bound is known, then the
corresponding lower bound is obtained as the previous floating-point
number, and vice versa. Thus, to implement interval semantics for
floating-point numbers is simpler than implementing interval
arithmetic in hardware.

But operands are not always point intervals: when 0 (alias $[0,m]$) or
$+\infty$ (alias $[M,+\infty)$) is an operand, the result may be a
wide interval. Then it is not the case that the upper bound is
unambiguously determined by the lower bound or vice versa. So for the
hardware implementation of floating-point arithmetic according to
interval semantics it is advantageous to be able to compute both
bounds as much as possible in parallel.

Many present day processors consist of multiple CPU cores along with a
separate floating point unit (FPU) that is IEEE 754 compliant
\cite{armfpu,osparc}. As with the \cite{osparc}, the CPU cores are
also multithreaded. With general-purpose processors being 754
compliant, all four of the rounding modes are implemented: round to
nearest, round to zero, round to $+\infty$ and round to $-\infty$. The
FPU is implemented using a pipeline architecture with the rounding
mode as part of the instruction word. This would allow the upper bound
and lower bound to be computed one cycle after each other. Thereby,
the computation of an interval result occurs within two cycles plus
the latency of the pipeline.  Before being fed to the FPU, the
appropriate operands are computed by the other CPU cores that exist
within the processor. The combination of these two operations yields
the correct arithmetic interval result.

To achieve interval semantics when the FPU only implements the
default rounding mode of round to nearest, a flag needs to be
available that indicates the direction of the rounding (there is no such
flag available in the processor's FPU).  In other words, after
rounding of the arithmetic result, the flag will show whether a round up
or truncation took place. In conjunction with the sign of the result,
correct rounding can be achieved in software or in the compiler. This
will require a small addition to the hardware, for which a flag can be
generated by comparing the last bit of the rounded word with the first
bit of the discarded portion of the pre-rounded word. If a pre-rounded
binary data of the mantissa is represented by
\begin{equation} 
   x =  b_{0}.b_{1}b_{2} \dots b_{r}b_{r+1} \dots b_{W} \\
\end{equation}
where $W$ is the pre-rounded word length and $r$ is the word length
of the data bus. Let $R\uparrow$ represent rounding up when set
active; the table below describes the result.
\vspace{.1in}
\begin{center}
\begin{tabular}{c | c || c }
	$b_{r}$ & $b_{r+1}$ & $R\uparrow$ \\ \hline
	0 & 0 & 0 \\
	0 & 1 & 1 \\
	1 & 0 & 0 \\
	1 & 1 & 1\\
\end{tabular}

\end{center}

\section{Conclusions}

We considered two possible interval semantics for the floating-point
zeroes: one with finite precision and one with infinite precision.
According to the first, there is no role for NaNs, which is
satisfying. Also, this semantics, possibly for the first time, gives a
mathematical definition for the distinction between the two
zeros. According to the infinite-precision semantics for the zeroes,
there is a role for NaNs, namely as name for the empty interval, and
this is also satisfying. Especially under this semantics it looks like
the framers of IEEE Standard 754 were guided by a keen intuition in
their choice of those outcomes where mathematics provided no guidance.
  
Numerical software should be the kind of software that is easiest to
verify, as it is modelled after numerical algorithms with proven
convergence properties. The IEEE standard accommodates overflow to a
certain extent by allowing the infinities as actual values. However,
it is difficult to verify the code in the presence of undefined values
for several operations on the infinities. Interval semantics gives a
mathematical interpretation to the standard floating-point values and
operations that is consistent with the mathematics of the reals, yet
avoids undefined results. Moreover, the distinction between $+0$ and
$-0$ has always been problematic from a mathematical point of
view. Interval semantics with finite-precision zeroes gives these
quantities distinct meanings that are consistent with the meanings of
the non-zero floating-point values.

\section{Acknowledgments}
We acknowledge support from the Natural Science and
Engineering Research Council NSERC.


\end{document}